\documentclass[letterpaper, 10 pt, conference]{ieeeconf}  

\IEEEoverridecommandlockouts                              
\usepackage{graphics} 
\usepackage{epsfig} 
\usepackage{times} 
\usepackage{amsmath} 
\usepackage{amssymb}  
\usepackage{algorithm}
\usepackage{algcompatible}

\usepackage{graphicx}
\usepackage{caption} 
\usepackage{subfig} 
\usepackage{url}

\usepackage{xcolor}
\definecolor{USred}{rgb}{0.74,0.1,0.1}
\definecolor{USblue}{rgb}{0.2,0.2,0.7}

\newcommand{\al}[1]{\begin{align} #1 \end{align}}

\newcommand{\G}{\mathcal{G}}

\newcommand{\Rs}{\mathbb{R}}
\newcommand{\Es}{\mathbb{E}}
\newcommand{\Gc}{\mathcal{G}}
\newcommand{\nn}{\nonumber}
\newcommand{\vct}{\mathrm{vec}}
\newcommand{\abs}{\mathrm{abs}}
\newcommand{\one}{\mathbf{1}}
\newcommand{\argmin}{\operatornamewithlimits{argmin}}

\DeclareMathOperator{\tr}{tr}

\newtheorem{problem}{Problem}

\title{\LARGE \bf Learning Quasi-Kronecker Product Graphical Models}

\author{Mattia Zorzi
\thanks{}
\thanks{M. Zorzi is with the Department of Information Engineering, University of Padova, Padova, Italy; email:	 
	 {\tt\small zorzimat@dei.unipd.it}}
\thanks{}%
}

\begin{document}

\maketitle
\thispagestyle{empty}
\pagestyle{empty}

\begin{abstract}
We consider the problem of learning graphical models where the support of the concentration matrix can be decomposed as a Kronecker product. We propose a method that uses the Bayesian hierarchical learning modeling approach. Thanks to the particular structure of the graph, we use a the number of hyperparameters which is small compared to the number of nodes in the graphical model. In this way, we avoid overfitting in the estimation of the hyperparameters. Finally, we test the effectiveness of the proposed method by a numerical example.
\end{abstract}

\section{INTRODUCTION}\label{sec:intro}

Many modern applications are characterized by high-dimensional data sets from which it is important to discover the meaningful interactions among the variables rather than to find an accurate model. An powerful tool to analyze these interrelations is given by graphical models (i.e. Markov networks), \cite{LAURITZEN_1996,speed1986gaussian,8378239}. The simplest version of the latter is constituted by a zero mean Gaussian random vector to which we attach an undirected graph: each node corresponds to a component of the random vector and there is an edge between two nodes if and only if the corresponding components are conditionally dependent given the others.

In these applications, there is a large interest of learning sparse graphical models (i.e. graphs with few edges) from data; indeed, these models are characterized by few  conditional interdependence relations among the components. 
Interestingly, a sparse graph corresponds to a covariance matrix whose inverse, say concentration matrix, is sparse. The problem of learning sparse graphical models, sometime called covariance selection problem, can be faced by using regularization techniques, \cite{huang2006covariance,banerjee2008model,friedman2008sparse,d2008first}. For instance, \cite{banerjee2008model} proposed a regularized maximum-likelihood (ML) estimator for the covariance matrix where the $\ell_1$ penalty norm on the concentration matrix has been considered. Since the $\ell_1$ norm penalty induces sparsity, the estimated covariance matrix will have a sparse inverse. It is worth noting that these approaches can be extended to dynamic graphical models, \cite{SONGSIRI_TOP_SEL_2010,LATENTG} as well as factor models \cite{valeCDC,ciccone2017factor,7331087,8264253}.

These regularized estimators are known to be sensitive to the choice of the regularization parameter, i.e. the weight on $\ell_1$ penalty, which is typically selected by cross-validation or theoretical derivation. To overcome this issue, a Bayesian hierarchical modeling approach has been considered, \cite{asadi2009map}. Here, the concentration matrix is modeled as a random matrix whose prior is characterized by a regularization parameter (called hyperparameter).  Then, the hyperparameter as well as the covariance matrix are jointly estimated. Since the $\ell_1$ norm shrinks all the entries to zero, and thus introduces a bias, a further improvement is to consider  
a weighted $\ell_1$ norm, see \cite{scheinberg2010sparse}, where the hyperparameter is a matrix whose dimension (in principle) coincides with the number of the nodes in the graph. On the other hand, the introduction of an hyperparameter with many variables could lead to overfitting in the estimation of the hyperparameter matrix.

An important class of graphical models is represented by the so called  Kronecker Product (KP) graphical models, \cite{tsiligkaridis2013covariance,4392825,dutilleul1999mle,tsiligkaridis2013convergence,SINQUIN201714131} wherein it is required that the concentration matrix can be decomposed as a Kronecker product. KP graphical models find application in many fields: spatiotemporal MEG/EEG modeling \cite{bijma2005spatiotemporal}; recommendation systems like NetFlix and gene expression analysis, \cite{allen2010transposable}; face recognition analysis\cite{zhang2010learning}. In these applications the most important feature is the graphical structure, i.e. the fact that the support of the concentration matrix can be decomposed as a Kronecker product.

The contribution of the present paper is to address the problem of learning graphical models where the support of the concentration matrix can be decomposed as a Kronecker product. We call such models Quasi-Kronecker Product (QKP) graphical models. Note that, the assumption that the support can be decomposed as a Kronecker product does not imply that the concentration matrix does. Therefore, QKP graphical models can understood as a weaker version of KP graphical models, making the former class less restrictive than the latter. Adopting the Bayesian hierarchical modeling approach, in the spirit of \cite{scheinberg2010sparse}, we  introduce two hyperparameter matrices whose total number of variables is small compared to the number of nodes in the graph. In this way, we avoid overfitting in the estimation of the hyperparameter.

The paper is outlined as follows. In Section \ref{sec:graph_model} we introduce graphical models and the problem of graphical model selection. In Section  \ref{sec:QK_graph} we introduce QKP graphical models. In Section  \ref{sec:learning} we propose a Bayesian procedure to learn QKP graphical models  from data, while Section \ref{sec:init} is devoted on how to initialize the procedure. In Section \ref{sec:sim} we present a numerical example to show the effectiveness of the proposed method. Finally, Section  \ref{sec:concl} draws the conclusions.

We warn the reader that the present paper only reports some preliminary result regarding the Bayesian estimation of QKP  graphical models. In particular, all the proofs and most of the technical assumptions needed therein are omitted and will be
published afterwards.

{\em Notation}: Given a symmetric matrix S, we write $S\succ 0$ ($S\succeq 0$) if $S$ is positive (semi-)definite. $x\sim \mathcal N(\mu,\Sigma)$ means $x$ is a Gaussian random vector with mean $\mu$ and covariance matrix $\Sigma$.
$\Es[\cdot]$ denotes the expectation operator. Given two functions $f(x)$ and $g(x)$, $f\propto g$ means that the argmin with respect to $x$ of $f$ and $g$ do coincide. Given a matrix $S$ of dimension $m\times m$, $s_{jk}$ denotes its entry in position $(j,k)$. Given a matrix $S$ of dimension $m_1m_2\times m_1m_2$,
$s_{jk,il}$ denotes its entry in position $((j-1)m_2+i,(k-1)m_2+l)$. Given a matrix $S$, $\vct(S)$ denotes the vectorization of matrix $S$. Given a matrix $S$ with positive entries, $log(S)$ denotes the matrix with entry $\log(s_{jk})$ in position $(j,k)$. Given a matrix $S$, $\exp(S)$ is the matrix with entry $\exp(s_{jk})$ in position $(j,k)$ and $\mathrm{abs}(S)$ denotes the matrix with entry $|s_{jk}|$ in position $(j,k)$. $\one_m$ denotes the $m$-dimensional vector of ones.

\section{GRAPHICAL MODEL SELECTION} \label{sec:graph_model}

\begin{figure}[htbp]
\centering
\includegraphics[width=0.5\columnwidth]{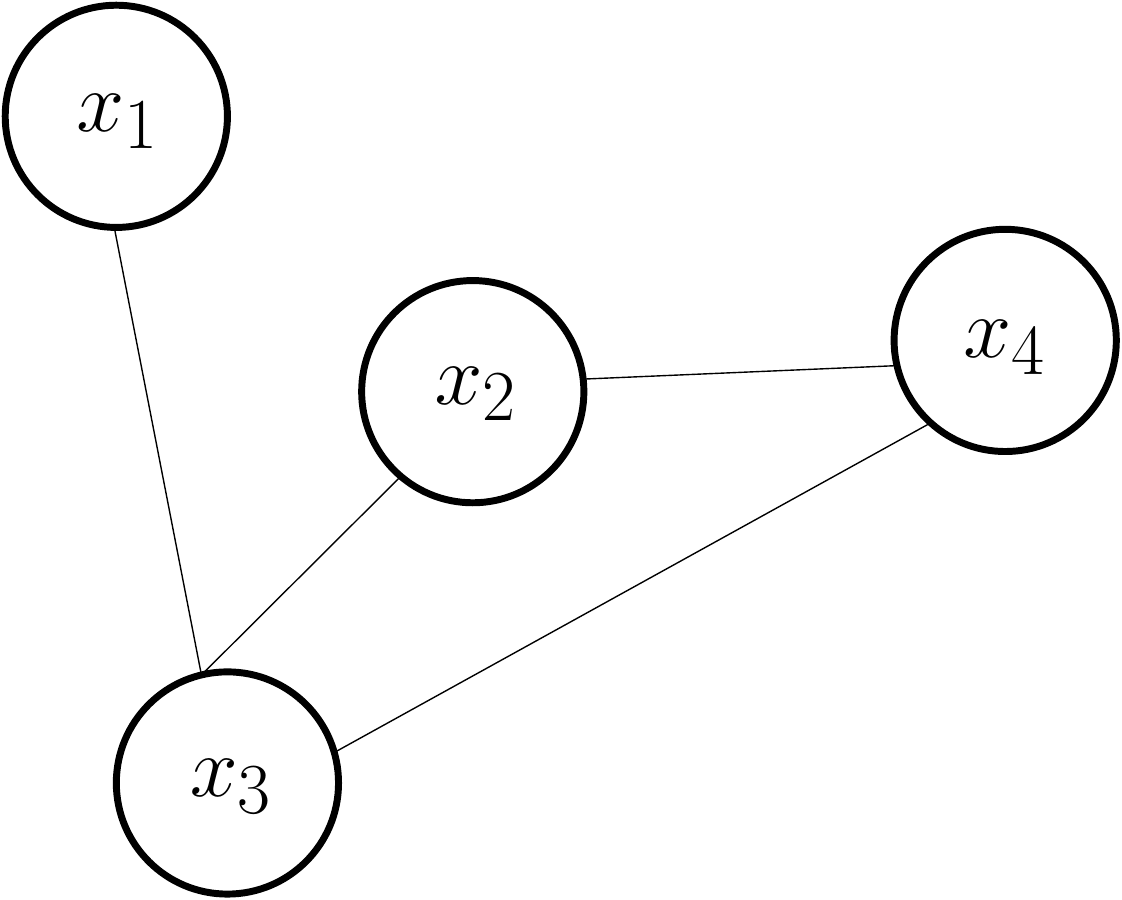}
\caption{Example of graphical model of $x=[\,x_1\,x_2\,x_3\,x_4\,]^T$ where $(x_1,x_3)$,
$(x_2,x_3)$, $(x_2,x_4)$ and  $(x_3,x_4)$ are the components of $x$ which are conditionally dependent given the others.}\label{fig_ex_grahp}
\end{figure}

Let $x=[\, x_1\ldots x_m\,]^T$ be a zero mean Gaussian random vector taking values in $\Rs^m$ and with covariance matrix $\Sigma\succ 0$. Thus, this random vector is completely characterized by $\Sigma$. We can attach to $x$ an undirected graph $\Gc(V,\Omega)$ where $V$ denotes the set of its nodes, and $\Omega$ denotes the set of its edges. More precisely,  each  
nodes corresponds to a component $x_j$, $j=1\ldots m$, of $x$ and there is an edge between nodes $j,k$ 
if and only $x_j$ and $x_k$ are conditionally dependent given the other components, or equivalently, for any $j\neq k$:
\al{x_j\, \bot 	,x_k\,|\,x_l, \, l\neq j,k\, \, \iff \, (j,k)\notin \Omega.}
Thus, $\Omega$ describes the conditionally dependent pairs of $x$. The graph $\Gc$ is referred to as graphical model of $x$; an example is provided in Figure 	\ref{fig_ex_grahp}.
Dempster proved that conditional independence relations are given by the concentration matrix of $x$, i.e. $S:=\Sigma^{-1}$, \cite{DEMPSTER_1972}:
\al{x_j\, \bot \, x_k\,|\,x_l, \, l\neq j,k\, \, \iff \, s_{jk}=0.}
Accordingly, sparsity of $S$, i.e. $S$ with many entries equal to zero, reflects the fact that the graphical model $\Gc$ of $x$ is sparse, i.e. $\Gc$ has few edges.

In many applications, it is required to learn a sparse graphical model $\Gc$ from data. More precisely, given a sequence of data $\mathrm x^N:=\{\,\mathrm x_1^T\ldots \mathrm x_N^T\,\}$ generated by $x$, find a sparse graphical model $\Gc(V,E)$ for $x$ where $\Gc$ is sparse. The simplest idea is to compute the sample covariance from the data 
\al{\hat \Sigma=\frac{1}{N}\sum_{k=1}^N \mathrm x_k \mathrm x_k^T;} then the graphical model is given by the support of $\hat \Sigma^{-1}$. However, the resulting graph  is full even in the case that the underlying system is well described by a sparse graphical model. In \cite{asadi2009map}, a procedure based on a Bayesian hierarchical model has been proposed. 
More precisely,  the entries of $S$ are assumed to be i.i.d. and Laplace distributed with hyperparameter $\gamma\geq 0$, i.e. the probability density function (pdf) of $s_{jk}$ is $p(s_{jk})=\gamma /2\exp\left(-\gamma|s_{jk}|\right)$. The resulting procedure is described in Algorithm \ref{algoS}. The main drawback of this approach is that it assigns a priori the same level of sparsity to each entry of $S$. This method has been extended to the case wherein only the entries in the same column of $S$ have the same distribution, \cite{scheinberg2010sparse}, allowing different levels of sparsity in the prior. A further extension is to assume that all the entries of $S$ may be distributed in a different way, but respecting the symmetry, that is $p(s_{jk})=p(s_{kj})=\gamma_{jk} /2\exp\left(-\gamma_{jk}|s_{jk}|\right)$ with $\gamma_{jk}\geq 0$. Using argumentations similar to the ones in \cite{scheinberg2010sparse}, it is not difficult to find that the procedure described in Algorithm \ref{algoS2}.

\begin{algorithm}
\caption{Algorithm for sparse graphical models in \cite{asadi2009map}}\label{algoS}
\begin{algorithmic}[1] \small
\STATE  Set $\varepsilon>0$, $\hat \gamma^{(0)}>0$, $\epsilon_{STOP}>0$
\STATE $h\leftarrow 1$
\REPEAT
\STATE Solve
\al{\footnotesize\hat S^{(h)}&=\underset{S\succ 0}{\, \mathrm{argmin}\,} -\frac{N}{2}\log |S| +\frac{N}{2} \tr(S\hat \Sigma)\nn\\ &\hspace{1.3cm}+\hat \gamma^{(h-1)}\sum_{j,k=1}^{m}| s_{jk}|\nn}
\STATE Update hyperparameter
\al{\hat \gamma^{(h)}&=\frac{m^2}{\sum_{j,k=1}^{m} |\hat s_{jk}^{(h)}| +\varepsilon}\nn}
\STATE $h\leftarrow h+1$
\UNTIL{ $\|\hat S^{(h)}-\hat S^{(h-1)}\|\leq \epsilon_{STOP}$}
\end{algorithmic}
\end{algorithm}

\begin{algorithm}
\caption{Algorithm for sparse graphical models}\label{algoS2}
\begin{algorithmic}[1] \small
\STATE  Set $\varepsilon>0$, $\hat \Gamma^{(0)}$ $m\times m$ matrix positive entrywise, $\epsilon_{STOP}>0$
\STATE $h\leftarrow 1$
\REPEAT
\STATE Solve
\al{\footnotesize\hat S^{(h)}&=\underset{S\succ 0}{\, \mathrm{argmin}\,} -\frac{N}{2}\log |S| +\frac{N}{2} \tr(S\hat \Sigma)\nn\\ &\hspace{1.3cm}+\sum_{j,k=1}^{m}\hat \gamma^{(h-1)}_{jk}| s_{jk}|\nn}
\STATE Update hyperparameter
\al{\hat \gamma_{jk}^{(h)}&=\frac{1}{|\hat s_{jk}^{(h)}| +\varepsilon}, \; j,k=1\ldots m\nn}
\STATE $h\leftarrow h+1$
\UNTIL{ $\|\hat S^{(h)}-\hat S^{(h-1)}\|\leq \epsilon_{STOP}$}
\end{algorithmic}
\end{algorithm}

\section{QKP GRAPHICAL MODELS}\label{sec:QK_graph}
\begin{figure*}[htbp]
\centering
\subfloat[]{ \includegraphics[width=0.5\columnwidth]{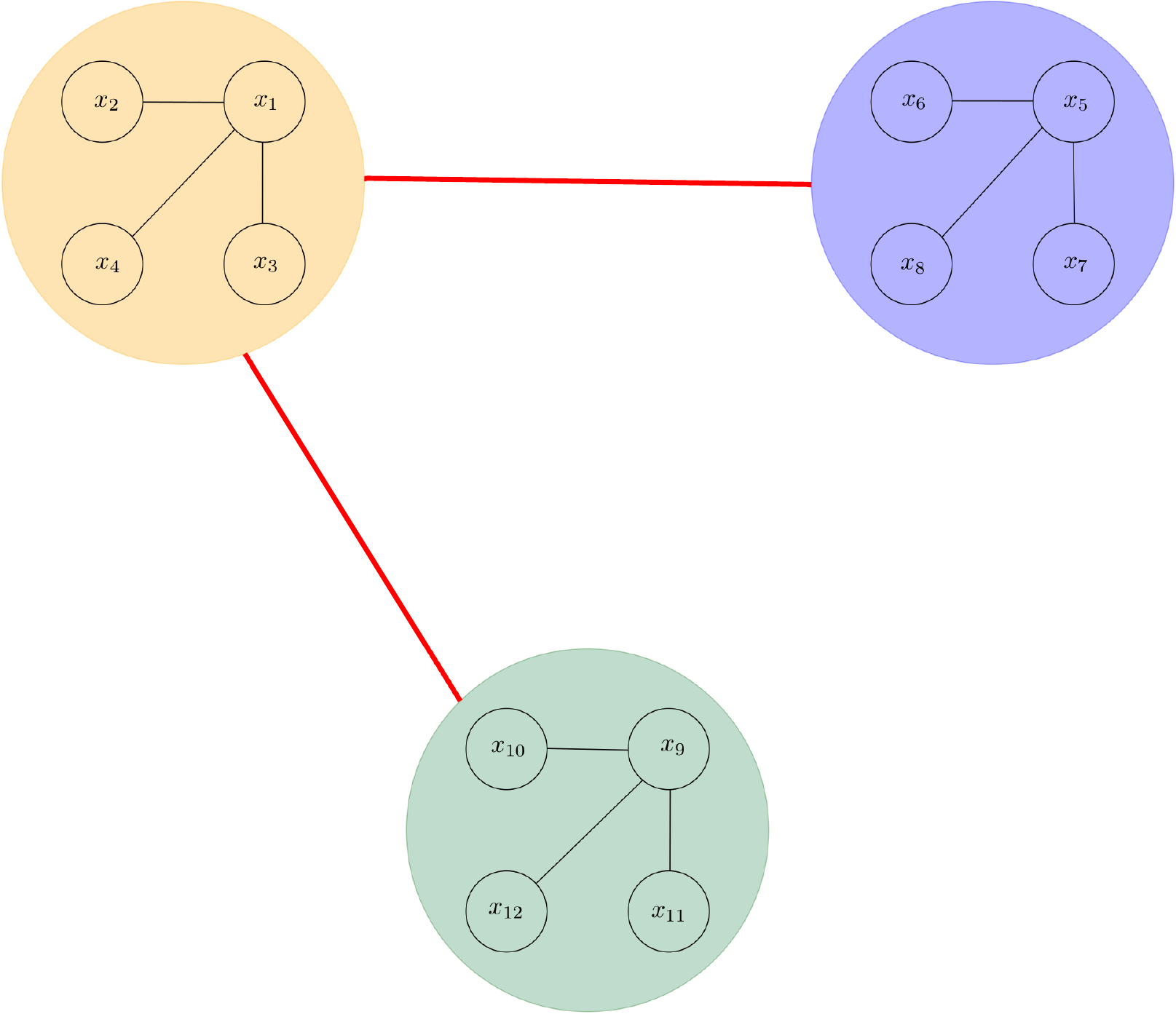} } \; \;
\subfloat[]{ \includegraphics[width=0.5\columnwidth]{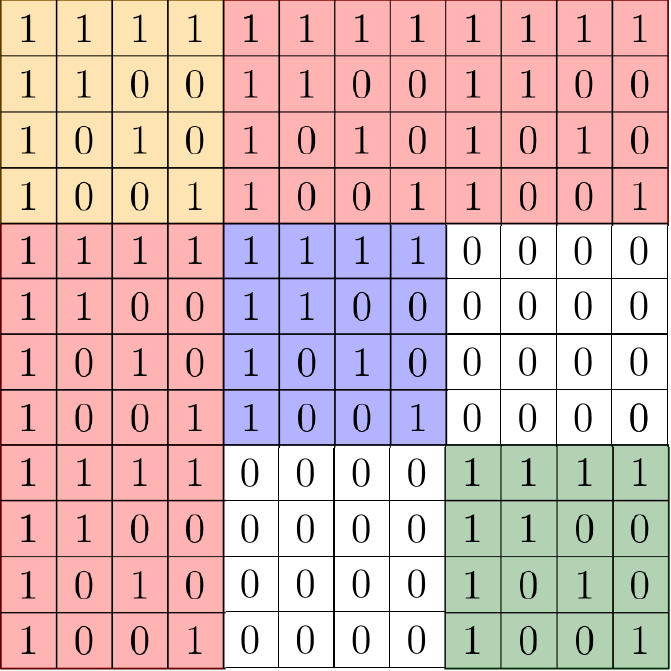}}\; \;
\subfloat[]{ \includegraphics[width=0.5\columnwidth]{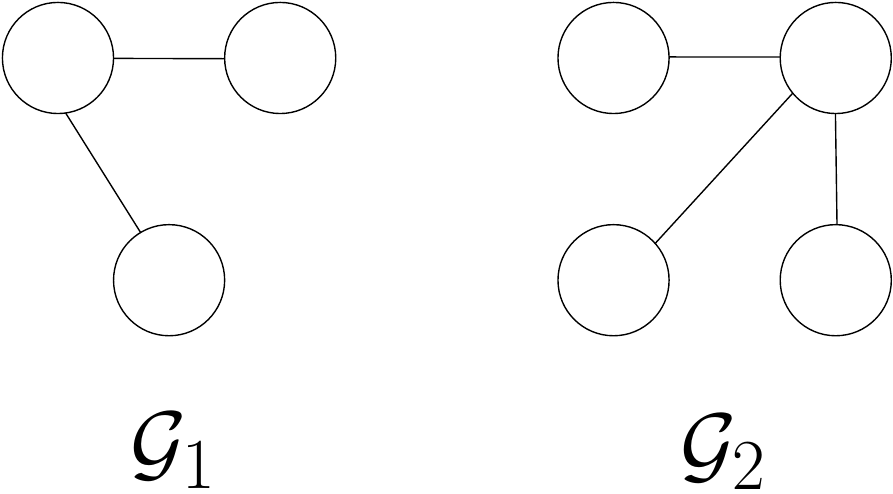}}
\caption{{\em Panel (a)}. Example of a Kronecker Product graph: we have three modules (yellow, blue, green): each module contains four nodes. Black edges describe conditional dependence relations among the nodes of the same module. Red edges describe conditional dependence relations among the modules.
{\em Panel (b)}. The corresponding matrix $E_\Omega$.
{\em Panel (c)}.  Graphs $\Gc_1(V_1,\Omega_1)$  and $\Gc_1(V_2,\Omega_2)$ generating the Kronecker Product graph $\Gc_1\otimes \Gc_2$ in Panel (a).} \label{Fig:ex_graph_kron}
\end{figure*}

Consider the undirected graph $\Gc(V,\Omega)$ and let $m$ denote the number of its nodes. Let $E_\Omega$ be the $m\times m$ binary matrix defined as follows:
\al{(E_\Omega)= \left\{\begin{array}{cc}1, & \hbox{ if $(j,k)\in\Omega$ }   \\ 0, &\hbox{ otherwise}  \end{array}\right. .}
We say that $\Gc(V,\Omega)$ is a Kronecker Product graph if there exist two graphs $\Gc_1(V_1,\Omega_1)$ and $\Gc_1(V_2,\Omega_2)$ with $m_1$ and $m_2$ nodes, respectively, such that:
\al{E_\Omega= E_{\Omega_1}\otimes E_{\Omega_2}}
where $m=m_1m_2$. In shorthand notation we will write $\Gc=\G_1\otimes \Gc_2$. In practice, in this graph $\Gc$ we can recognize modules containing $m_2$ nodes sharing the same graphical structure, described by $\Omega_2$; the interaction among those $m_1$ modules is described by $\Omega_1$. An illustrative example is given in Figure \ref{Fig:ex_graph_kron}.

Let $x=[\, x_1\ldots x_m\,]^T$ be a zero mean Gaussian random vector taking values in $\Rs^m$ and with inverse covariance matrix (i.e. concentration matrix) $S\succ 0$ whose support is $E_\Omega=E_{\Omega_1}\otimes E_{\Omega_2}$. We can attach to $x$ a Kronecker Product  graph $\Gc(V,\Omega)=\Gc_1(V_1,\Omega_1)\otimes \Gc_2(V_2,\Omega_2)$ Accordingly, $\Omega_1$ characterizes the conditional dependence relations among the modules, while $\Omega_2$ characterizes the recurrent conditional dependence relations among the nodes in each  module. This graphical model is referred to as Quasi-Kronecker Product (QKP) graphical model to distinguish from the Kronecker Product (KP) graphical model proposed in \cite{tsiligkaridis2013convergence}. Indeed, in the latter the support of the concentration matrix and its support admit a Kronecker decomposition, i.e.  $E_\Omega=E_{\Omega_1}\otimes E_{\Omega_2}$ and $S=S_1\otimes S_2$. In our graphical model, even if the
support of the concentration matrix admits a Kronecker product decomposition, the concentration matrix does not.

\section{LEARNING QKP GRAPHICAL MODELS} \label{sec:learning}

We address the problem of learning a QKP graphical model $\Gc=\Gc_1\otimes \Gc_2$ from data. In many real applications the observed data are explained by a sparse graphical model because the latter allows a straightforward interpretation of the interaction among the variables involved in the application. Thus, in our case we require that $\Gc_1$ and $\Gc_2$ are sparse. The problem can be formalized as follows.

\begin{problem}\label{probl_kron}Let $x$ be a zero mean random vector of dimension $m=m_1m_2$, $m_1,m_2\in\mathbb{N}$, with zero mean and inverse covariance matrix $S$. Given a sequence of data $\mathrm x^N:=\{\,\mathrm x_1^T\ldots \mathrm x_N^T\,\}$ generated by $x$, find a QKP graphical model $\Gc(V,E)=\Gc_1(V_1)\otimes \Gc_2(V_2,E_2)$ for $x$ where $\Gc_1$ and $\Gc_2$ are sparse and have $m_1,m_2$ nodes, respectively.\end{problem}

To solve Problem \ref{probl_kron}, we adopt the Bayesian hierarchical modeling. $S$ is modeled as a random matrix whose prior depends on the hyperparameters $\Lambda$ and $\Gamma$. $\Lambda$ is a $m_1\times m_1$ symmetric random matrix with nonnegative entries and $\Gamma$ is a $m_2\times m_2$ symmetric random matrix with nonnegative entries. The hyperprior for $\Lambda$ and $\Gamma$ depends on $\varepsilon_1$ and $\varepsilon_2$, respectively. The latter are deterministic positive quantities.

We proceed to characterize the Bayesian model in detail.
The conditional pdf of $x^N$ under model $x\sim \mathcal N(0,S^{-1})$ is:
\al{\label{pdf_Fish}p(\mathrm x^N|S)&=\prod_{k=1}^N\frac{1}{\sqrt{(2\pi)^m |S^{-1}|}} \exp\left(-\frac{1}{2}\mathrm x_k^T S\mathrm x_k^T \right). \nn\\
&\propto |S|^{N/2}\exp\left(-\frac{1}{2}\tr(S\hat \Sigma)\right) }
where the neglected terms do not depend on $S$. In what follows we assume that $\hat \Sigma\succ 0$. We model the entries of $S$ as independent random variables, so that the prior for $S$ is:
\al{p(S|\Lambda,\Gamma)=\prod_{j,k=1}^{m_1}\prod_{i,l=1}^{m_2} p(s_{jk,il}|\lambda_{jk},\gamma_{il})} 
where $s_{jk,il}$ is Laplace distributed
\al{\label{prima_oss}p(s_{jk,il}|\lambda_{jk},\gamma_{il})=\frac{\lambda_{jk}\gamma_{il}}{2} \exp\left(-\lambda_{jk}\gamma_{il}|s_{jk,il}|\right).} 
$\Lambda$ and $\Gamma$ are independent random matrices, i.e. $p(\Lambda,\Gamma)=p(\Lambda)p(\Gamma)$. The entries of $\Lambda$ and $\Gamma$ are assumed to be independent 
\al{p(\Lambda)=\prod_{j,k=1}^{m_1} p(\lambda_{jk}), \; \; p(\Gamma)=\prod_{i,l=1}^{m_2} p(\gamma_{il})}
and with exponential distribution
\al{\label{seconda_oss} p(\lambda_{jk})=\varepsilon_1 \exp(-\varepsilon_1 \lambda_{jk}),\;\; p(\gamma_{il})=\varepsilon_2 \exp(-\varepsilon_2 \gamma_{il})} with $\varepsilon_1,\varepsilon_2$ deterministic and positive quantities.  At this point, some comment regarding the choice of the prior on $S$ and the hyperprior on $\Lambda,\Gamma$ is required.  From (\ref{prima_oss}) it is clear that $s_{jk,il}$ takes value close to zero with high probability if the product $\lambda_{jk}\gamma_{il}$ is large. Moreover, if $\gamma_{il}$ is very large for some $(i,l)$ and $\lambda_{jk}\geq \epsilon>0$, for all $j,k=1\ldots m_1$, such that $\gamma_{il}\epsilon$ is large then $s_{jk,il}$ with  $j,k=1\ldots m_1$ take values close to zero with high probability. 
Accordingly, the different modules in the graph will have a similar sparsity pattern with high probability. Accordingly, prior (\ref{prima_oss}) assigns high probability to QKP graphical models. The hyperprior in (\ref{seconda_oss}) guarantees that $\lambda_{jk}$ and $\gamma_{il}$ diverge with probability zero. As we will see, this assumption guarantees that the optimization procedure that we propose is well-posed.

Next, we characterize the maximum a posteriori (MAP) estimator of $S$ (and thus also the MAP estimator of the covariance matrix by the invariance principle). The latter minimizes the negative log-likelihood 
\al{\ell (\mathrm x^N; S,\Lambda,\Gamma)= -\log p(\mathrm x^N,S, \Lambda,\Gamma)}
where $p(\mathrm x^N,S, \Lambda,\Gamma)$ is the joint pdf of $\mathrm x^N$, $S$, $\Lambda$ and $\Gamma$.
Note that,
\al{\label{joint_pdf}p(\mathrm x^N,S, \Lambda,\Gamma)=p(\mathrm x^N|S)p(S|\Lambda,\Gamma)p(\Lambda)p(\Gamma),} 
in particular the negative log-likelihood contains the prior (\ref{prima_oss}) inducing sparsity on intra-group/modules. Moreover, we have \al{\label{log_joint} \ell( &\mathrm x^N; S, \Lambda,\Gamma):=-\log p(\mathrm x^N|S)-\log p(S|\Lambda,\Gamma)\nn\\
& \hspace{0.5cm}-\log p(\Lambda)-\log p(\Gamma) \nn\\
 &\propto -\frac{N}{2}\log |S| +\frac{N}{2} \tr(S\hat \Sigma)+\sum_{j,k=1}^{m_1}\sum_{i,l=1}^{m_2} \lambda_{jk}\gamma_{il}|s_{jk,il}|\nn\\ 
 &\hspace{0.5cm} +\sum_{j,k=1}^{m_1}( \varepsilon_1 \lambda_{jk} -m_2^2\log\lambda_{jk})+\sum_{il=1}^{m_2}(\varepsilon_2\gamma_{il}- m_1^2\log\gamma_{il})
  }
 where the neglected terms do not depend on $S$, $\Lambda$ and $\Gamma$. It is clear that 
the MAP estimator of $S$ depends on $\Lambda$, $\Gamma$, $\varepsilon_1$ and $\varepsilon_2$. In what follows we assume $\varepsilon_1$ and $\varepsilon_2$ fixed. Then,  a way to estimate 
$\Lambda$ and $\Gamma$ from the data is provided by the empirical Bayes approach: $\Lambda$ and $\Gamma$ are given by maximizing the marginal likelihood
of $x^{N}$ which is obtained by integrating out $S$ in (\ref{joint_pdf}), \cite{rasmussen2004gaussian}. However, it is not easy to find an analytical expression for the marginal likelihood in this case.  
An alternative simplified approach for optimizing $\Lambda$, $\Gamma$ is the generalized maximum likelihood (GML) method, \cite{zhou1997approximate}. According to this method, $S$, $\Lambda$ and $\Gamma$ minimize jointly  (\ref{log_joint}):
\al{(\hat S,\hat \Lambda,\hat \Gamma)=&\underset{S,\Lambda,\Gamma}{\, \mathrm{argmin}\,} \ell( \mathrm x^N; S, \Lambda,\Gamma)\nn\\
& \hbox{ s.t. } S\succ 0,\; \lambda_{jk}\geq 0,\; \gamma_{il}\geq 0.} Since the joint optimization of the three variables is still an hard problem, we propose an iterative three-step procedure. At the $h$-th  iteration we solve the following three optimization problems:
 \al{\label{P1}\hat S^{(h)}&=\underset{S\succ 0}{\, \mathrm{argmin}\,} \ell( \mathrm x^N; S,\hat  \Lambda^{(h-1)},\hat \Gamma^{(h-1)})\\
\label{P2} \hat \Lambda^{(h)}&=\underset{\substack{\Lambda\\ \hbox{\footnotesize s.t. } \lambda_{jk}\geq 0}}{\, \mathrm{argmin}\,} \ell( \mathrm x^N; \hat S^{(h)}, \Lambda,\hat \Gamma^{(h-1)})\\
\label{P3}\hat \Gamma^{(h)}&=\underset{\substack{\Gamma\\ \hbox{\footnotesize s.t. } \gamma_{il}\geq 0}}{\, \mathrm{argmin}\,} \ell( \mathrm x^N; \hat S^{(h)}, \hat \Lambda^{(h)},\Gamma).}

It is possible to prove that Problems (\ref{P1}), (\ref{P2}) and (\ref{P3}) do admit a unique solution. Moreover, the resulting procedure is illustrated in Algorithm \ref{algo}.
\begin{algorithm}
\caption{Proposed algorithm}\label{algo}
\label{algo:RWS}
\begin{algorithmic}[1] \small
\STATE  Set $\varepsilon_1>0$, $\varepsilon_2>0$, $\hat \Lambda^{(0)}$ and $\hat \Gamma^{(0)}$ positive entrywise, $\epsilon_{STOP}>0$
\STATE $h\leftarrow 1$
\REPEAT
\STATE Solve
\al{\footnotesize\hat S^{(h)}&=\underset{S\succ 0}{\, \mathrm{argmin}\,} -\frac{N}{2}\log |S| +\frac{N}{2} \tr(S\hat \Sigma)\nn\\ &\hspace{1.3cm}+\sum_{j,k=1}^{m_1}\sum_{i,l=1}^{m_2} \hat \lambda_{jk}^{(h-1)}\hat \gamma_{il}^{(h-1)}|s_{jk,il}|\nn}
\STATE Update hyperparameters
\al{\hat \lambda_{jk}^{(h)}&=\frac{m_2^2}{\sum_{i,l=1}^{m_2} \hat\gamma_{il}^{(h-1)}|\hat s_{jk,il}^{(h)}| +\varepsilon_1}, \; \; j,k=1\ldots m_1\nn\\
\hat \gamma_{il}^{(h)}&=\frac{m_1^2}{\sum_{j,k=1}^{m_1} \hat \lambda_{jk}^{(h)}|\hat s_{jk,il}^{(h)}| +\varepsilon_2},\; \; i,l=1\ldots m_2\nn}
\STATE $h\leftarrow h+1$
\UNTIL{ $\|\hat S^{(h)}-\hat S^{(h-1)}\|\leq \epsilon_{STOP}$}
\end{algorithmic}
\end{algorithm}
In the aforementioned algorithm the hyperparameters selection is performed iteratively through Step 5. It is worth noting that Algorithm \ref{algo} is similar to Algorithm \ref{algoS} and Algorithm \ref{algoS2}: the main difference is that in the proposed algorithm we have two types of hyperparameters. As a consequence, in the proposed algorithm we have three optimization steps instead of two.

\section{INITIAL CONDITIONS} \label{sec:init}
In Algorithm \ref{algo} we have to fix the initial conditions for the hyperparameters, that is $\hat \Lambda^{(0)}$ and $\hat \Gamma^{(0)}$. 
The idea is to approximate $\hat \Sigma^{-1}$ through a Kronecker product, then the two matrices of this product are used to initialize $\Lambda$ and $\Gamma$.

Given $\hat \Sigma$, we want  to find $\bar W$ and $\bar Y$ of dimension $m_1$ and $m_2$, respectively, with positive entries such that $\bar W\otimes \bar Y \approx  \abs(\hat \Sigma^{-1})+\epsilon  \one_{m_1m_2} \one_{m_1m_2}^T $ where $\epsilon >0$ is chosen sufficiently small. The presence of the term $\epsilon \one_{m_1m_2} \one_{m_1m_2}^T$ allows to take the entrywise logarithm on both sides, obtaining 
\al{\label{approx1}	 W &\otimes I_{m_2} +I_{m_1}\otimes Y\nn\\ &\approx  \log (\abs(\hat \Sigma^{-1})+\epsilon  \one_{m_1m_2} \one_{m_1m_2}^T)} where $ W=\log \bar W$ and $ Y=\log \bar Y$.
It is not difficult to see that the following relations hold:
\al{\vct ( W \otimes I_{m_2})&=(I_{m_1}\otimes \one_{m_2}\otimes I_{m_1} \otimes \one_{m2})\vct(	W)\nn\\
\vct (I_{m_1} \otimes Y)&=( \one_{m_1}\otimes I_{m_2} \otimes \one_{m1}\otimes I_{m_2}) \vct( Y).} Then, we can write (\ref{approx1}) as $Az\approx b$ where 
\al{A&=[\, I_{m_1}\otimes \one_{m_2}\otimes I_{m_1} \otimes \one_{m2}\; \; \one_{m_1}\otimes I_{m_2} \otimes \one_{m1}\otimes I_{m_2}\,]	\nn\\
z&=[\,\vct( W)^T \;\; \vct( Y)^T\,]^T\nn\\
b&=\vct(\log (\abs(\hat \Sigma^{-1})+\epsilon  \one_{m_1m_2} \one_{m_1m_2}^T)).} Accordingly, $z$ can be found by solving the least squares problem $\hat z=\argmin_z \|Az-b\|$, therefore $\hat z=(A^TA)^{-1}A^Tb$. From $z$ we recover $W$ and $Y$. Note that, $W$ and $Y$ computed from $z$ are not symmetric matrices. Thus, to compute $\bar W$ and $\bar Y$ from $W$ and $Y$ we force the symmetric structure: $\bar W=(\exp (W)+\exp (W)^T)/2$ and $\bar Y =(\exp (Y)+\exp (Y)^T)/2$. At this point, it is worth noting that $w_{jk}y_{il}$ provides roughly the order of magnitude of $s_{jk,il}$. On the other hand, the hyperparameters $\lambda_{jk}$ and $\gamma_{il}$ provides the prior about the order the order of magnitude of $s_{jk,il}$: the larger $\lambda_{jk}\gamma_{il}$ is, the more $s_{jk,il}$ is close to zero. Accordingly, we choose $\hat \Lambda^{(0)}$ and $\hat \Gamma^{(0)}$ such that:
\al{\hat \lambda_{jk}^{(0)}&=\frac{1}{w_{jk}},\; j,k=1\ldots m_1\nn\\
\hat \gamma_{il}^{(0)}&=\frac{1}{y_{il}},\; i,l=1\ldots m_2 .}

\section{SIMULATION RESULTS} \label{sec:sim}
 \begin{figure*}[htbp]
\centering
\includegraphics[width=2\columnwidth]{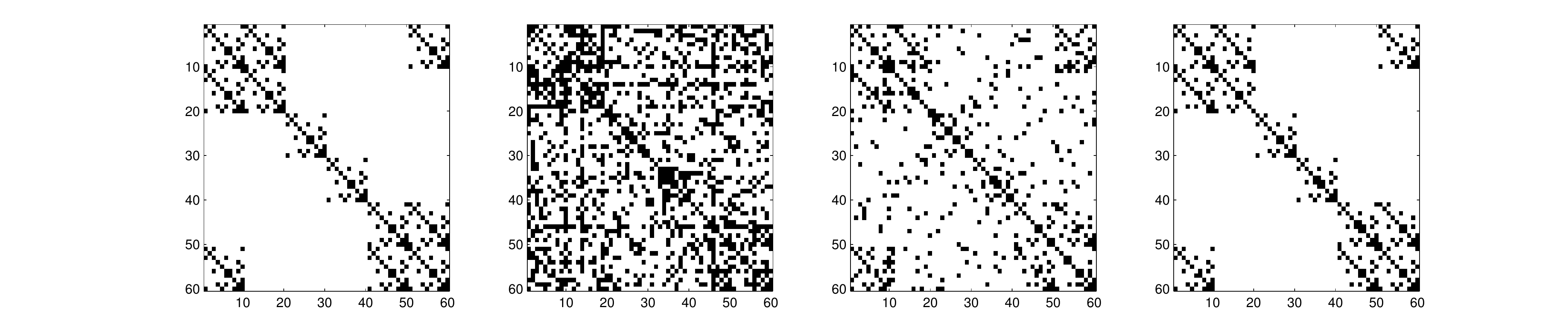}
\caption{Set of edges estimated in a realization of the Monte Carlo experiment. {\em First picture}. $E_\Omega$ where $\Omega$ is the set of edges of the true model. {\em Second picture}. $E_{\hat \Omega_{S1}}$ where $\hat \Omega_{S1}$ is the set of edges estimated using S1. {\em Third picture}. $E_{\hat \Omega_{S2}}$ where $\hat \Omega_{S2}$ is the set of edges estimated using S2.  {\em Fourth picture}. $E_{\hat \Omega_{QKP}}$ where $\hat \Omega_{QKP}$ is the set of edges estimated using QKP. }\label{fig:realization1000}
\end{figure*} 
We consider a Monte Carlo experiment structured as follows:
\begin{itemize}
\item We generate $60$ QKP graphical models with $m_1=6$ modules and each module contains $m_2=10$ nodes. For each model, $\Omega_1$ and $\Omega_2$ are generated randomly. The fraction of edges is set equal to $20\%$ for both $\Omega_1$ and $\Omega_2$;
\item For each model we generate a finite-length realization $\mathrm x^N:=\{\,\mathrm x_1\ldots \mathrm x_N\,\}$, with $N=1000$, and we compute the sample covariance $\hat \Sigma$.   
\item For each realization we consider the following estimators:  
\begin{itemize}
\item S1 estimator: it computes a sparse graphical model by using Algorithm  \ref{algoS}, in this case we have one scalar hyperparameter;
\item S2 estimator: it computes a sparse graphical model by using Algorithm  \ref{algoS2}, in this case we have $m_1m_2(m_1m_2+1)/2=1830$ variables in the  hyperparameter matrix;
\item QKP estimator: it computes a QKP graphical model with $m_1=6$ modules and each module has $m_2=10$ nodes, in this case the total number of variables in the hyperparameters matrices is  $m_1(m_1+1)/2+m_1(m_2+1)/2=76$.
\end{itemize}
\item For each realization,  we compute the relative error in reconstructing the concentration matrix  and the relative error in reconstructing the sparsity pattern using S1, S2 and QKP. For instance, the relative error in reconstructing concentration matrix using QKP estimator is 
\al{e_{QKP}=\frac{\| S_{true}- \hat S_{QKP}\|}{\|S_{true}\|}}
where $S_{true}$ denotes the concentration matrix of the true model, while $\hat S_{QKP}$ is the estimated concentration matrix; here, $\|\cdot \|$ denotes the Frobenius norm. The relative error in reconstructing the sparsity pattern using QKP estimator is 
\al{e_{SP,QK}=\frac{\| E_{\Omega_{true}}- E_{\hat \Omega_{QKP}}\|}{(m_1m_2(m_1m_2+1))/2}}
where $\Omega_{set}$ and $\hat \Omega_{QKP}$ denote the set of edges of the true graphical model and the one estimated, respectively.
\end{itemize} 

\begin{figure}[htbp]
\centering
\includegraphics[width=\columnwidth]{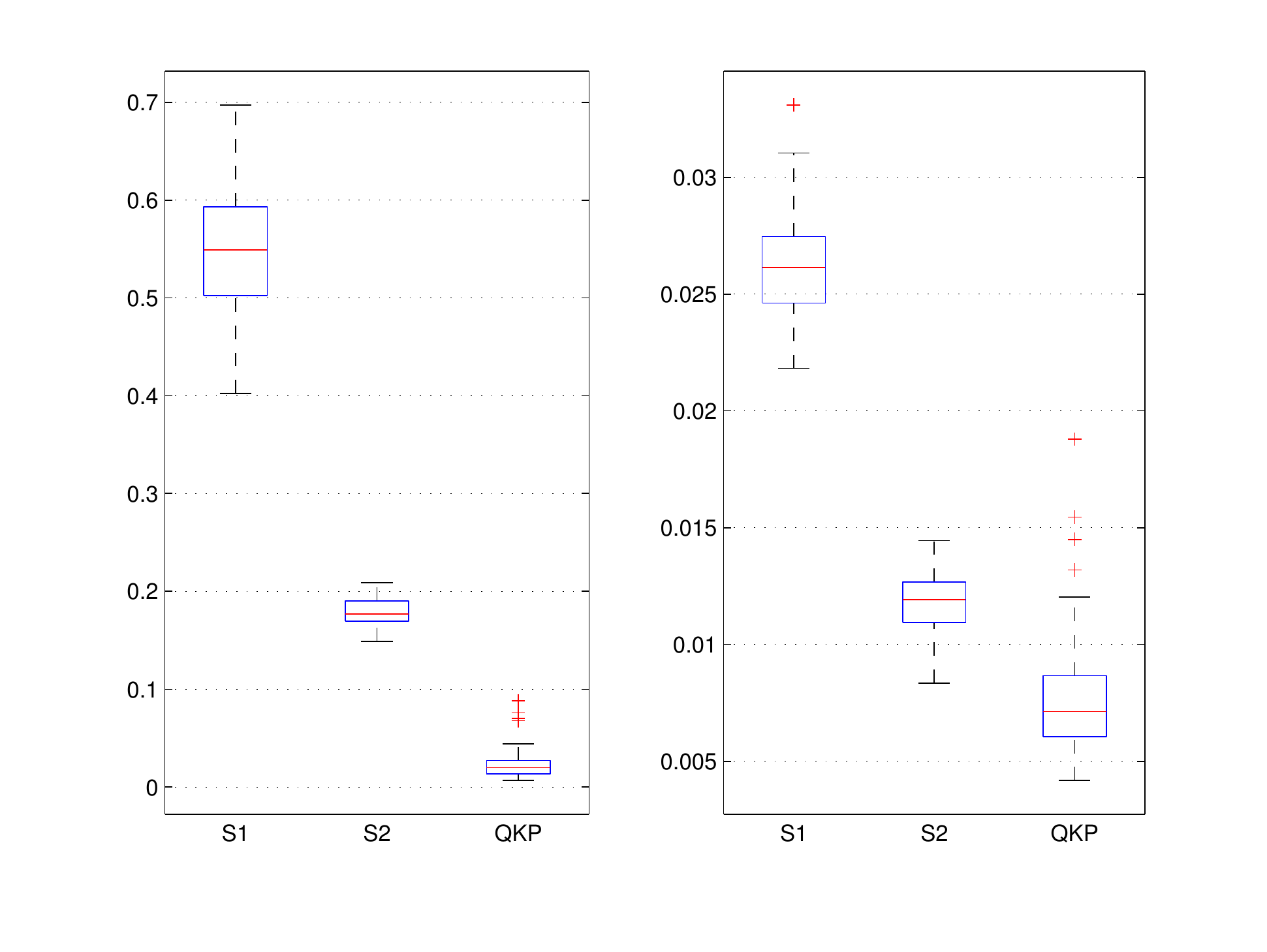}
\caption{Monte Carlo experiment with realizations of length $N=1000$. {\em Left panel}. Boxplot of the relative error in reconstructing the sparsity pattern. {\em Right panel}. Boxplot of  the relative error in reconstructing the concentration matrix. }\label{fig:boxplot1000}
\end{figure}

Figure \ref{fig:realization1000} depicts the set of edges estimated using the three estimators in a realization of the Monte Carlo experiment. As we see, only QK provides a
structure similar to the true model. Figure \ref{fig:boxplot1000}  shows the boxplot of the relative error in reconstructing the sparsity pattern (left panel) and the concentration matrix  (right panel).
As we can see, the worst performance is given by S1, while the best performance is achieved by QKP. In particular, the relative error of the estimated sparsity pattern for QKP is very  small compared to the other two methods. The poor performance of S1 is due by the fact that only a scalar hyperparameter is not sufficient to capture the correct structure of the graph. On the contrary, the poor performance of S2 is due by the fact that there is overfitting in the estimation of the hyperparameter matrix.

\section{CONCLUSIONS} \label{sec:concl}
We have introduced Quasi-Kronecker Product graphical models wherein the nodes are regrouped in modules having the same number of nodes. The interactions among the nodes of the same module as well as the interaction among the nodes of two modules follow a common structure. Then, we have addressed the problem of learning QKP graph models from data using a Bayesian hierarchical model. Finally, we have compared the proposed procedure with Bayesian learning techniques for estimating sparse graphical models: simulation evidence showed the effectiveness of the proposed method.

\end{document}